\documentclass[a4paper,11pt]{article}
                                                                                                                                                                                                                                                                                                                                                                                                                                                                                                                                                           
\usepackage{amssymb,amsthm,amscd, amsbsy, array,amsmath,mathdots}
\usepackage[francais]{babel}
\usepackage{fontenc}
\usepackage{amsfonts}
\usepackage[mathscr]{eucal}
\usepackage{mathrsfs}
\usepackage{latexsym}
\usepackage{graphicx}
\usepackage[a4paper,textwidth=17cm,textheight=21cm]{geometry}

\newcommand{\be}{\begin{equation}}
\newcommand{\ee}{\end{equation}}
\newcommand{\bea}{\begin{eqnarray}}
\newcommand{\eea}{\end{eqnarray}}

\def\1{{\mbox{\boldmath $1$}}} %
\newcommand{\C}{\mathbb{C}}
\begin{document}
\thispagestyle{empty}
\begin{center}
	\begin{huge}
		Some Poisson-Lie sigma models
	\end{huge}
\end{center}

\begin{center}
	\begin{large}
		ARM Boris \\
	\end{large}
	{\footnotesize  \textit{arm.boris@gmail.com}}
\end{center}

\vspace{5cm}
\begin{center}
	\textbf{Abstract} \\
	We calculate the Poisson-Lie sigma model for every $4$-dimensional\\ Manin triples (function of its structure constant) and we
	\\ give the $6$-dimensional models for the Manin triples\\ $(\mathfrak{sl}(2,\C)\oplus \mathfrak{sl}(2,\C)^*,\mathfrak{sl}(2,\C),\mathfrak{sl}(2,\C)^*)$,\\
	$(\mathfrak{sl}(2,\C)\oplus \mathfrak{sl}(2,\C)^*,\mathfrak{sl}(2,\C)^*),\mathfrak{sl}(2,\C)$, \\
	$(\mathfrak{sl}(2,\C),\mathfrak{su}(2,\C),\mathfrak{sb}(2,\C))$ and\\
	  $(\mathfrak{sl}(2,\C),\mathfrak{sb}(2,\C),\mathfrak{su}(2,\C))$
\end{center}

\newpage

\setlength{\textheight}{21cm}
\addtolength{\voffset}{0.5cm}

\section{Introduction}
A Manin triples $(\mathfrak{D},\mathfrak{g},\tilde{\mathfrak{g}})$ is a bialgebra $(\mathfrak{g},\tilde{\mathfrak{g}}$ which don't intersect each others and a direct sum of this bialgebra $\mathfrak{D}=\mathfrak{g}\oplus \tilde{\mathfrak{g}})$. If the corresponding Lie groups have a Poisson structure, they are called Poisson-Lie groups.
A Poisson-Lie sigma models is an action (\ref{actiondeux}) calculated by a Poisson vector field matrix.
\cite{vysoky} have deduced the extremal field which minimize the action of this models, which gives the motion equation (\ref{equmotion}).
We calculate here the action and the equations of motion for some $6$-dimensionals Manin triples and we give  a general formula for each $4$-dimensional Manin triples.
The $6$-dimensional Manin triples are
 $(\mathfrak{sl}(2,\C)\oplus \mathfrak{sl}(2,\C)^*,\mathfrak{sl}(2,\C),\mathfrak{sl}(2,\C)^*)$,$(\mathfrak{sl}(2,\C)\oplus \mathfrak{sl}(2,\C)^*,\mathfrak{sl}(2,\C)^*),\mathfrak{sl}(2,\C)$,$(\mathfrak{sl}(2,\C),\mathfrak{su}(2,\C),\mathfrak{sb}(2,\C))$ and $(\mathfrak{sl}(2,\C),\mathfrak{sb}(2,\C),\mathfrak{su}(2,\C))$.

\newpage
\section{Some Manin triples}
The Drinfeld double $D$ is defined as a Lie group such that its Lie algebra $\mathfrak{D}$ equipped by a symmetric ad-invariant nondegenerate bilinear form $<.,.>$ can be decomposed into a pair of maximally isotropic subalgebras  $\mathfrak{g}, \tilde{\mathfrak{g}}$ such that $\mathfrak{D}$  as a vector space is the direct sum of $\mathfrak{g}$ and $\tilde{\mathfrak{g}}$. Any such decomposition written as an ordered set $(\mathfrak{D},\mathfrak{g},\tilde{\mathfrak{g}})$ is called a Manin triples $(\mathfrak{D},\mathfrak{g},\tilde{\mathfrak{g}})$,$(\mathfrak{D},\tilde{\mathfrak{g}},\mathfrak{g})$. 
\paragraph{}
One can see that the dimensions of the subalgebras are equal and that bases $\{ T_i\} ,\{\tilde{T}^i\}$ in the subalgebras can be chosen so that 
\begin{equation}
<T_i,T_j>=0, \hspace{3mm}<T_i,\tilde{T}^j>=<\tilde{T}^j,T_i>=\delta_i^j, \hspace{3mm}<\tilde{T}^i,\tilde{T}^j>=0
\end{equation}
This canonical form of the bracket is invariant with respect to the transformations 
\begin{equation}
T_i'=T_k A^k_i, \hspace{3mm} \tilde{T}'^j=(A^{-1})_k^j \tilde{T}^k
\end{equation}
Due to the ad-invariance of $<.,.>$ the algebraic structure of $\mathfrak{D}$ is 
\begin{eqnarray}
\nonumber [T_i,T_j]=c_{ij}^{\hspace{3mm}k} T_k, &&\hspace{3mm}[\tilde{T}^i,\tilde{T}^j]=f^{ij}_{\hspace{3mm}k} \tilde{T}^k \\
\nonumber [T_i, \tilde{T}^j]&=&f^{jk}_{\hspace{3mm}i} T_k -c_{ik}^{\hspace{3mm}j} \tilde{T}^k
\end{eqnarray}
\paragraph{}
 There are just  four types of nonisomorphic four-dimensional Manin triples.
\\
\textit{Abelian Manin triples:}
\begin{equation}
[T_i,T_j]=0, \hspace{6mm} [\tilde{T}^i, \tilde{T}^j]=0, \hspace{6mm} [T_i, \tilde{T}^j]=0,\hspace{6mm} i,j=1,2
\end{equation}
\textit{Semi-Abelian Manin triples (only non trivial brackets are displayed):}
\begin{equation}
 [\tilde{T}^1, \tilde{T}^2]=\tilde{T}^2, \hspace{6mm} [T_2, \tilde{T}^1]=T_2,\hspace{6mm} [T_2, \tilde{T}^2]=-T_1
\end{equation}
\textit{Type A non-Abelian Manin triples ($\beta \neq 0$):}
\begin{eqnarray}
\nonumber  [T_1,T_2]=T_2,&& \hspace{6mm}[\tilde{T}^1, \tilde{T}^2]=\beta \tilde{T}^2\\
\nonumber   [T_1, \tilde{T}^2]=-\tilde{T}^2,\hspace{6mm} [T_2, \tilde{T}^1]&&=\beta T_2\hspace{6mm},[T_2, \tilde{T}^2]=-\beta T_1+\tilde{T}^1
\end{eqnarray}
\textit{Type B non-Abelian Manin triples:}
\begin{eqnarray}
\nonumber  [T_1,T_2]=T_2, &&\hspace{6mm}[\tilde{T}^1, \tilde{T}^2]= \tilde{T}^1\\
\nonumber  [T_1, \tilde{T}^1]=T_2,\hspace{6mm} [T_1, \tilde{T}^2]&=&
  -T_1-\tilde{T}^2\hspace{6mm},[T_2, \tilde{T}^2]=\tilde{T}^1
\end{eqnarray}
\paragraph{}
Now we focus some six dimensional Manin triples. We recall that the commutation relations of the Lie algebra $\mathfrak{sl}(2,\C)$ of the Lie group $SL(2,\C)$:
\begin{equation}
[ T_1,T_2]=2T_2 , \hspace{3mm}[T_1,T_3]=-2T_3,\hspace{3mm}[T_2,T_3]=T_1
\end{equation}
The dual Lie algebra $\mathfrak{sl}(2,\C)^*$ of the Lie algebra $\mathfrak{sl}(2,C)$ has the commutation relations:
\begin{equation}
[ \tilde{T}^1,\tilde{T}^2]=\frac{1}{4}\tilde{T}^2 ,\hspace{3mm}[ \tilde{T}^1,\tilde{T}^3]=\frac{1}{4}\tilde{T}^3, \hspace{3mm}[ \tilde{T}^2,\tilde{T}^3]=0
\label{relcommslet}
\end{equation}
There is a scalar product on $(\mathfrak{sl}(2,\C)\oplus \mathfrak{sl}(2,\C)^*$ such that  (see \cite{kosmann}):
\begin{equation}
(T_i,\tilde{T}^j)=\delta_i^j
\end{equation}
Finally, we have that $(\mathfrak{sl}(2,\C)\oplus \mathfrak{sl}(2,\C)^*,\mathfrak{sl}(2,\C),\mathfrak{sl}(2,\C)^*)$ with this scalar product  is a Manin triple.
We note that $(\mathfrak{sl}(2,\C)\oplus \mathfrak{sl}(2,\C)^*,\mathfrak{sl}(2,\C)^*,\mathfrak{sl}(2,\C))$ with this scalar product is also a Manin triples.
\paragraph{}
The Iwasawa decomposition allows us to decompose:
\begin{equation}
\mathfrak{sl}(2,\C)=\mathfrak{su}(2,\C)\oplus \mathfrak{sb}(2,\C)
\label{Iwasawa}
\end{equation}
where $\mathfrak{su}(2,\C)$ is the Lie algebra of the Lie group $SU(2)$  with commutation relations:
\begin{equation}
[ T_1,T_2]=T_3 , \hspace{3mm}[T_2,T_3]=T_1,\hspace{3mm}[T_3,T_1]=T_2
\end{equation}
$\mathfrak{sb}(2,\C)$ is the Lie algebra of the Borel subgroup $SB(2,\C)$ with commutation relations:
\begin{equation}
[ \tilde{T}^1,\tilde{T}^2]=\tilde{T}^2 ,\hspace{3mm}[ \tilde{T}^1,\tilde{T}^3]=\tilde{T}^3, \hspace{3mm}[ \tilde{T}^2,\tilde{T}^3]=0
\label{relcommsb}
\end{equation}
Here we can see in comparing (\ref{relcommsb}) and (\ref{relcommslet}) that $\mathfrak{sb}(2,\C)\simeq \mathfrak{sl}(2,\C)^*$.\\

The Iwasawa decomposition (\ref{Iwasawa}) allows us to identify $\mathfrak{sb}(2,\C)\simeq \mathfrak{su}(2,\C)^*$.
We define a scalar product on $(\mathfrak{sl}(2,\C)\oplus \mathfrak{sl}(2,\C)^*$ such that $(x,y)=$Im$($Tr$(x | y))$.
With this scalar product we have (see \cite{kosmann}):
\begin{equation}
(T_i,\tilde{T}^j)=\delta_i^j
\end{equation}
Finally we have that 
$(\mathfrak{sl}(2,\C),\mathfrak{su}(2,\C),\mathfrak{sb}(2,\C))$ with this scalar product  is a Manin triple.
We note that $(\mathfrak{sl}(2,\C),\mathfrak{sb}(2,\C),\mathfrak{su}(2,\C))$ with this scalar product is also a Manin triples.

\section{Poisson-Lie sigma models}
Given a Lie group $M$ and a Poisson structure on it. We define the action of this model (see \cite{hajizadeh}) as:
\begin{equation}
S_1=\int_\Sigma (< dg g^{-1} ,A>-\frac{1}{2}<A,(r-Ad_g r Ad_{g})A>)
\end{equation}
where $g \in G,A=A^i_\alpha d\xi^\alpha X_i$ and $r \in \mathfrak{g}\otimes 
\mathfrak{g}$ is a classical $r$ matrix with $\mathfrak{g}$ as the Lie algebra of $G$ and $\{ X_i\}$ as a basis of $\mathfrak{g}$. Note that the above action can be applied for simple or nonsemisimple Lie group $G$ with ad-invariant symmetric bilinear nondegenerate form $<X_i,X_j>=G_{ij}$ on the Lie algebra $\mathfrak{g}$. When the metric $G_{ij}$ of Lie algebra is denegerate then the above action is not good. Here we use the following action instead of the above one:
\begin{equation}
S_2= \int_\Sigma (dX_i \wedge A_i -\frac{1}{2}\mathcal{P}^{ij} A_i \wedge A_j)
\label{actiondeux}
\end{equation}
where $xî$ are Lie group parameters with parametrization  (e.g.)
\begin{equation}
\forall g \in G , g=e^{X_1 T_1}e^{X_2 T_2}...
\end{equation}
where $P^{ij}$ is the Poisson structure on the Lie group which for coboundary Poisson Lie groups it is obtained from
\begin{equation}
(\mathcal{P}(g))_\chi =b(g)a(g)^{-1}
\label{matricedepoisson}
\end{equation}
We can obtain $a(g)^{-1}$ and $b(g)$ in computing:
\begin{eqnarray}
(Ad_{g^{-1}})_\chi &=&
\left(\begin{array}{cc}
a(g)^T& b(g)^T\\
0&d(g)^T\\   
   \end{array}
\right)
\\
(Ad_g)_\chi &=&
\left(\begin{array}{cc}
a(g)^{-T}& -a(g)^{-T}b(g)^T d(g)^{-T}\\
0&d(g)^{-T}\\   
   \end{array}
\right)
\label{matriceactionadjointe}
\end{eqnarray}
where $^T$ denotes the transpose matrix.

The extremal fiels $(X,A)$ which minimize the action
(\ref{actiondeux}) have to satisfy the equation written locally (see \cite{vysoky}) as:
\begin{eqnarray}
dX_i+\mathcal{P}^{ij}(X)A_j=0\\
dA_k+\frac{1}{2}\mathcal{P}^{ij}_{\hspace{3mm},k}(X) A_i \wedge A_j =0
\label{equmotion}
\end{eqnarray}
where $\mathcal{P}^{ij}_{\hspace{3mm},k}=\partial_k \mathcal{P}^{ij}|_{X_k=0}$.

\section{Poisson-Lie sigma model of any $4$-dimensional Manin triple}
We first calculate the matrix of the adjoint actions function of structure constant:
\begin{eqnarray}
\nonumber  ad_{T_1}&=& 
\left(\begin{array}{cccc}
0&c_{12}^{\hspace{3mm}1}&0&-f^{12}_{\hspace{3mm}1}\\
0&c_{12}^{\hspace{3mm}2}&f^{12}_{\hspace{3mm}1}&0\\   
0&0&0&0\\
0&0&-c_{12}^{\hspace{3mm}1}&-c_{12}^{\hspace{3mm}2}\\   
   \end{array}
\right)
\\
\nonumber ad_{T_2}&=& 
\left(\begin{array}{cccc}
-c_{12}^{\hspace{3mm}1}&0&0&-f^{12}_{\hspace{3mm}2}\\
-c_{12}^{\hspace{3mm}2}&0&f^{12}_{\hspace{3mm}2}&0\\   
0&0&c_{12}^{\hspace{3mm}1}&c_{12}^{\hspace{3mm}2}\\
0&0&0&0\\   
   \end{array}
\right)
\end{eqnarray}
To obtain the matrix $\mathcal{P}$, we calculate the adjoint action matrix of a general element $g=\prod_{i=1}^2 e^{\alpha_i T_i}$ by the formula:
\begin{equation}
(Ad_{\prod_{i=1}^2 e^{X_i T_i}})_\chi =\prod_{i=1}^2 e^{X_i (ad_{T_i})_\chi}
\end{equation}
Similarily, we have:
\begin{equation}
(Ad_{(\prod_{i=1}^2 e^{X_i T_i})^{-1}})_\chi =\prod_{i=1}^2 e^{-X_{3-i} (ad_{T_{3-i}})_\chi}
\end{equation}
We can deduce 
the matrix $\mathcal{P}^{ij}$:
\begin{eqnarray}
\nonumber  \mathcal{P}^{ij}&=& 
\left(\begin{array}{ccc}
0&-\mathcal{P}^{21}\\
\mathcal{P}^{21}&0\\   
   \end{array}
\right)
\end{eqnarray}
where
\begin{equation}
\mathcal{P}^{21}=\frac{c_{12}^{\hspace{3mm}1}(-1+e^{c_{12}^{\hspace{3mm}2} X_1})f^{12}_{\hspace{3mm}1}+c_{12}^{\hspace{3mm}2}e^{c_{12}^{\hspace{3mm}2} X_1-c_{12}^{\hspace{3mm}1} X_2}(-1+e^{c_{12}^{\hspace{3mm}1} X_2})f^{12}_{\hspace{3mm}2}}{c_{12}^{\hspace{3mm}2}c_{12}^{\hspace{3mm}1}}
\label{generalformula}
\end{equation}
Now, we can calculate the action (\ref{actiondeux}) of the model 
\begin{equation}
S_2 = \int_\Sigma \sum_{i=1}^2 dX_i \wedge A_i -\mathcal{P}^{21}A_2 \wedge A_1
\end{equation}
and the equations of motion (\ref{equmotion}):
\begin{eqnarray}
\nonumber dX_1-\mathcal{P}^{21}A_2=0 \\
\nonumber dX_2+\mathcal{P}^{21}A_1=0\\
\nonumber dA_1-\frac{c_{12}^{\hspace{3mm}1}c_{12}^{\hspace{3mm}2} f^{12}_{\hspace{3mm}1}+c_{12}^{\hspace{3mm}2}c_{12}^{\hspace{3mm}2}(-e^{-c_{12}^{\hspace{3mm}1} X_2}+1)f^{12}_{\hspace{3mm}2}}{c_{12}^{\hspace{3mm}2}c_{12}^{\hspace{3mm}1}} A_1 \wedge A_2  =0 \\
\nonumber dA_2+\frac{c_{12}^{\hspace{3mm}2}e^{c_{12}^{\hspace{3mm}2} X_1}c_{12}^{\hspace{3mm}1}f^{12}_{\hspace{3mm}2}}{c_{12}^{\hspace{3mm}2}c_{12}^{\hspace{3mm}1}}A_2 \wedge A_1=0\\
\end{eqnarray}

\section{Poisson-Lie sigma model of $(\mathfrak{sl}(2,\C)\oplus \mathfrak{sl}(2,\C)^*,\mathfrak{sl}(2,\C),\mathfrak{sl}(2,\C)^*)$}
We first calculate the matrix of the adjoint actions:
\begin{eqnarray}
\nonumber  ad_{T_1}&=& 
\left(\begin{array}{cccccc}
0&0&0&0&0&0\\
0&2&0&0&0&0\\   
0&0&-2&0&0&0\\
0&0&0&0&0&0\\   
0&0&0&0&-2&0\\
0&0&0&0&0&2\\   
   \end{array}
\right)
\\
\nonumber ad_{T_2}&=\frac{1}{4}& 
\left(\begin{array}{cccccc}
0&0&4&0&-1&0\\
-8&0&0&1&0&0\\   
0&0&0&0&0&0\\
0&0&0&0&8&0\\   
0&0&0&0&0&0\\
0&0&0&-4&0&0\\   
   \end{array}
\right)
\\
\nonumber ad_{T_3}&=\frac{1}{4}& 
\left(\begin{array}{cccccc}
0&-4&0&0&0&-1\\
0&0&0&0&0&0\\   
8&0&0&1&0&0\\
0&0&0&0&0&-8\\   
0&0&0&4&0&0\\
0&0&0&0&0&0\\   
   \end{array}
\right)
\end{eqnarray}
To obtain the matrix $\mathcal{P}$, we calculate the adjoint action matrix of a general element $g=\prod_{i=1}^3 e^{\alpha_i T_i}$ by the formula:
\begin{equation}
(Ad_{\prod_{i=1}^3 e^{X_i T_i}})_\chi =\prod_{i=1}^3 e^{X_i (ad_{T_i})_\chi}
\end{equation}
Similarily, we have:
\begin{equation}
(Ad_{(\prod_{i=1}^3 e^{X_i T_i})^{-1}})_\chi =\prod_{i=1}^3 e^{-X_{4-i} (ad_{T_{4-i}})_\chi}
\end{equation}
We can deduce 
the matrix $\mathcal{P}^{ij}$:
\begin{eqnarray}
\nonumber  \mathcal{P}^{ij}&=& 
\left(\begin{array}{ccc}
0&-\frac{X_2}{4}(1+X_2 X_3)e^{2X_1}&-\frac{X_3}{4}e^{-2X_1}\\
\frac{X_2}{4}(1+X_2 X_3)e^{2X_1}&0&\frac{X_2X_3}{2}\\   
\frac{X_3}{4}e^{-2X_1}&-\frac{X_2X_3}{2}&0\\
   \end{array}
\right)
\end{eqnarray}
Now, we can calculate the action (\ref{actiondeux}) of the model 
\begin{equation}
S_2 = \int_\Sigma \sum_{i=1}^3 dX_i \wedge A_i +(\frac{X_2}{4}(1+X_2 X_3)e^{2X_1})A_1 \wedge A_2 +\frac{X_3}{4}e^{-2X_1} A_1 \wedge A_3 -\frac{X_2X_3}{2} A_2 \wedge A_3
\end{equation}
and the equations of motion (\ref{equmotion}):
\begin{eqnarray}
\nonumber dX_1-(\frac{X_2}{4}(1+X_2 X_3)e^{2X_1})A_2-\frac{X_3}{4}e^{-2X_1} A_3=0 \\
\nonumber dX_2+(\frac{X_2}{4}(1+X_2 X_3)e^{2X_1})A_1+\frac{X_2X_3}{2}  A_3=0\\
\nonumber dX_3+\frac{X_3}{4}e^{-2X_1} A_1-\frac{X_2X_3}{2}  A_2=0\\
\nonumber dA_1 -\frac{X_2}{2}(1+X_2 X_3)A_1 \wedge A_2 +\frac{X_3}{2} A_1\wedge A_2 =0 \\
\nonumber dA_2- \frac{e^{2X_1}}{4} A_1 \wedge A_2 +\frac{X_3}{2}A_2\wedge A_3=0\\
\nonumber dA_3-\frac{X_2^2}{4}e^{2X_1} A_1 \wedge A_2 -\frac{e^{-2 X_1}}{4}A_1 \wedge A_3 +X_2 A_2 \wedge A_3=0
\end{eqnarray}

\section{Poisson-Lie sigma model of $(\mathfrak{sl}(2,\C)\oplus \mathfrak{sl}(2,\C)^*,\mathfrak{sl}(2,\C)^*,\mathfrak{sl}(2,\C))$}
Now to obtain this Poisson Lie sigma model, we have to change 
$T_i \to \tilde{T}^i$ and $\tilde{T}^i \to T_i$ of the previous model.
And we can calculate  
the matrix of the adjoint actions as we do previously.
With this we can deduce 
the matrix $\mathcal{P}^{ij}$ for this model:
\begin{eqnarray}
\nonumber  \mathcal{P}^{ij}&=& 
\left(\begin{array}{ccc}
0&-2e^{\frac{X_1}{4}}X_2&2e^{\frac{X_1}{4}}X_3\\
2e^{\frac{X_1}{4}}X_2&0&2-\frac{1}{2}e^{\frac{X_1}{2}}(4+X_2X_3)\\   
-2e^{\frac{X_1}{4}}X_3&-2+\frac{1}{2}e^{\frac{X_1}{2}}(4+X_2X_3)&0\\
   \end{array}
\right)
\end{eqnarray}
Now, we can calculate the action (\ref{actiondeux}) of the model 
\begin{equation}
S_2 = \int_\Sigma \sum_{i=1}^3 dX_i \wedge A_i +2 e^{\frac{X_1}{4}}X_2 A_1 \wedge A_2 -2 e^{\frac{X_1}{4}}X_3 A_1 \wedge A_3 +(-2+\frac{1}{2}e^{\frac{X_1}{2}}(4+X_2X_3)) A_2 \wedge A_3
\end{equation}
and the equations of motion (\ref{equmotion}):
\begin{eqnarray}
\nonumber dX_1-2e^{\frac{X_1}{4}}X_2 A_2+2e^{\frac{X_1}{4}}X_3  A_3=0 \\
\nonumber dX_2+2e^{\frac{X_1}{4}}X_2 A_1+(2-\frac{1}{2}e^{\frac{X_1}{2}}(4+X_2X_3))  A_3=0\\
\nonumber dX_3-2e^{\frac{X_1}{4}}X_3 A_1-(2-\frac{1}{2}e^{\frac{X_1}{2}}(4+X_2X_3))  A_2=0\\
\nonumber dA_1 -\frac{X_2}{2}A_1 \wedge A_2 +\frac{X_3}{4} A_1\wedge A_3-\frac{1}{4}(4+X_2X3)A_2 \wedge A_3 =0 \\
\nonumber dA_2-2 e^{\frac{X_1}{4}} A_1 \wedge A_2 -\frac{1}{2} e^{\frac{X_1}{4}} A_2\wedge A_3=0\\
\nonumber dA_3+2 e^{\frac{X_1}{4}} A_1 \wedge A_3 -\frac{e^{\frac{X_1}{2}}X_2}{2} A_2 \wedge A_3 =0
\end{eqnarray}

\section{Poisson-Lie sigma model of $(\mathfrak{sl}(2,\C),\mathfrak{su}(2,\C),\mathfrak{sb}(2,\C))$}
We can calculate  
the matrix of the adjoint actions as we do previously.
With this we can deduce 
the matrix $\mathcal{P}^{ij}$ for this model:
\begin{eqnarray}
\nonumber  \mathcal{P}^{ij}&=& 
\left(\begin{array}{ccc}
 \scriptstyle 0&  \scriptstyle  -\cos X_1 \cos X_3 \sin X_2 +\sin X_1 \sin X_3 & \scriptstyle -\cos X_3 \sin X_1 \sin X_2-\cos X_1 \sin X_3 \\
 \scriptstyle \cos X_1 \cos X_3 \sin X_2-\sin X_1 \sin X_3 & 0&  \scriptstyle -1+\cos X_2 \cos X_3\\   
 \scriptstyle \cos X_3 \sin X_1 \sin X_3+\cos X_1 \sin X_3& \scriptstyle 1-\cos X_2 \cos X_3& \scriptstyle 0\\
   \end{array}
\right)
\end{eqnarray}
Now, we can calculate the action (\ref{actiondeux}) of the model 
\begin{eqnarray}
\nonumber S_2 = \int_\Sigma \sum_{i=1}^3 dX_i \wedge A_i -(-\cos X_1 \cos X_3 \sin X_2 +\sin X_1 \sin X_3 )A_1 \wedge A_2 \\
-(-\cos X_3 \sin X_1 \sin X_2-\cos X_1 \sin X_3) A_1 \wedge A_3 -( -1+\cos X_2 \cos X_3)  A_2 \wedge A_3
\end{eqnarray}
and the equations of motion (\ref{equmotion}):
\begin{eqnarray}
\nonumber dX_1+(-\cos X_1 \cos X_3 \sin X_2 +\sin X_1 \sin X_3 ) A_2+(-\cos X_3 \sin X_1 \sin X_2-\cos X_1 \sin X_3)  A_3=0 \\
\nonumber dX_2+(\cos X_1 \cos X_3 \sin X_2-\sin X_1 \sin X_3) A_1+(-1+\cos X_2 \cos X_3)  A_3=0\\
\nonumber dX_3+(\cos X_3 \sin X_1 \sin X_3+\cos X_1 \sin X_3) A_1(1-\cos X_2 \cos X_3) A_2=0\\
\nonumber dA_1 +\sin  X_3 A_1\wedge A_2 -\cos X_3 \sin X_2 A_1 \wedge A_3=0 \\
\nonumber dA_2-\cos X_1 \cos X_3 A_1\wedge A_2-\cos X_3 \sin X_1 A_1 \wedge A_3=0\\
\nonumber dA_3 +\sin X_1 A_1 \wedge A_2 -\cos X_1 A_1 \wedge A_3 =0
\end{eqnarray}

\section{Poisson-Lie sigma model of $(\mathfrak{sl}(2,\C),\mathfrak{sb}(2,\C),\mathfrak{su}(2,\C))$}
We can calculate  
the matrix of the adjoint actions as we do previously.
With this we can deduce 
the matrix $\mathcal{P}^{ij}$ for this model:
\begin{eqnarray}
\nonumber  \mathcal{P}^{ij}&=& 
\left(\begin{array}{ccc}
0&-e^{ X_1}X_3 & -e^{ X_1} X_2\\
e^{ X_1} X_3 & 0&\frac{1}{2}(1-e^{2 X_1}(1+X_2^2+X_3^2))\\   
e^{X_1} X_2 &\frac{1}{2}(-1+e^{2X_1}(1+X_2^2+X_3^2))&0\\
   \end{array}
\right)
\end{eqnarray}
Now, we can calculate the action (\ref{actiondeux}) of the model 
\begin{equation}
S_2 = \int_\Sigma \sum_{i=1}^3 dX_i \wedge A_i +e^{X_1} X_3 A_1 \wedge A_2+e^{X_1}X_2  A_1 \wedge A_3 -\frac{1}{2}(1-e^{2X_1}(1+2X_2^2+2X_3^2))  A_2 \wedge A_3
\end{equation}
and the equations of motion (\ref{equmotion}):
\begin{eqnarray}
\nonumber dX_1-e^{X_1}X_3 A_2-e^{X_1}X_2  A_3=0 \\
\nonumber dX_2+e^{X_1}X_3 A_1+\frac{1}{2}(1-e^{2X_1}(1+X_2^2+X_3^2))  A_3=0\\
\nonumber dX_3+e^{X_1}X_2 A_1-\frac{1}{2}(1-e^{2X_1}(1+X_2^2+X_3^2))  A_2=0\\
\nonumber dA_1 -X_3 A_1\wedge A_2 -X_2 A_1 \wedge A_3 -(1+X_2^2+X_3^2)A_2 \wedge A_3=0 \\
\nonumber dA_2- e^{X_1} A_1\wedge A_3=0\\
\nonumber dA_3 -e^{X_1} A_1 \wedge A_2=0
\end{eqnarray}

\section{Discussion}
We gives here the Poisson-Lie sigma models of some Manin triples. Concerning the general formula (\ref{generalformula}), we have to say that this is no problem when $c_{12}^{\hspace{3mm}2}$ and $c_{12}^{\hspace{3mm}2}$ is zero because
\begin{equation}
\mathcal{P}^{21}=\frac{(-1+e^{c_{12}^{\hspace{3mm}2} X_1})f^{12}_{\hspace{3mm}1}}{c_{12}^{\hspace{3mm}2}}+\frac{e^{c_{12}^{\hspace{3mm}2} X_1-c_{12}^{\hspace{3mm}1} X_2}(-1+e^{c_{12}^{\hspace{3mm}1} X_2})f^{12}_{\hspace{3mm}2}}{c_{12}^{\hspace{3mm}2}}
\label{generalformula}
\end{equation}
which can be approximate by
\begin{equation}
\mathcal{P}^{21}=(X_1+\frac{c_{12}^{\hspace{3mm}2}}{2}X_1^2+...)f^{12}_{\hspace{3mm}1}+e^{c_{12}^{\hspace{3mm}2} X_1-c_{12}^{\hspace{3mm}1} X_2}(X_2+\frac{c_{12}^{\hspace{3mm}1}}{2}X_2^2+...)f^{12}_{\hspace{3mm}2}
\label{generalformula}
\end{equation}
We tried to obtain the equivalent formula for $n=3$ but the calculus was too hard.


\newpage
\addcontentsline{toc}{chapter}{Bibliographie}
\begin{thebibliography}{9}
\bibitem{hajizadeh}
Hajizadeh S., Rezaei-Aghdam A., Poisson-Lie Sigma models over low dimensional real Poisson-Lie groups
\bibitem{kosmann}
Kosmann-Schwarzbach Y., Lie bialgebras, Poisson Lie Groups and Dressing Transformation
\bibitem{vysoky}
Vysok\'y J., Hlavat\'y , Poisson Lie Sigma Models on Drinfeld double
\end{thebibliography}
\end{document}